\documentclass[12pt,a4paper]{article}
\usepackage{a4,amssymb,graphics,epsfig,latexsym}
\usepackage[centertags]{amsmath}
\usepackage{amsfonts}
\def\be{\begin{equation}}
\def\bea{\begin{eqnarray*}}
\def\ee{\end{equation}}
\def\eea{\end{eqnarray*}}
\def\ba{\begin{array}}
\def\ea{\end{array}}
\def\bi{\begin{itemize}}
\def\ei{\end{itemize}}
\newtheorem{theo}{Theorem}

\newcommand{\Lip}{{\rm Lip}}
\newcommand{\SU}{{\rm SU}}
\newcommand{\su}{{\rm su}}
\def\GG#1{\mbox{\gothic #1}}
\font\gothic=eufm10 scaled 1100

\def\ZZ{{{\rm Z}\kern-.4em{\rm Z}}}
\def\RR{{{\rm I}\kern-.2em{\rm R}}}
\def\NN{{{\rm I}\kern-.2em{\rm N}}}
\def\TT{{{\rm T}\kern-.5em{\rm T}}}
\def\CC{{{\rm I}\kern-.5em{\rm C}}}

\begin{document}
\title{Close-to-Optimal Bounds for
$\mathbf{\SU}(N)$ Loop Approximation}
\author{Peter Oswald\\
{\small Jacobs University Bremen}\\
{\small e-mail: p.oswald@jacobs-university.de} \and
Tatiana Shingel\\
{\small DAMTP, Cambridge University}\\{\small e-mail:
t.shingel@damtp.cam.ac.uk}}

\date{}
\maketitle

\begin{abstract}
In \cite{OS}, we proved an asymptotic O$(n^{-\alpha/(\alpha+1)})$
bound for the approximation of $\SU(N)$ loops ($N\ge 2$) with
Lipschitz smoothness $\alpha>1/2$ by polynomial loops of degree
$\le n$. The proof combined factorizations of $\SU(N)$ loops into
products of constant $\SU(N)$ matrices and loops of the form
$e^{A(t)}$ where $A(t)$ are essentially $\su(2)$ loops
preserving the Lipschitz smoothness, and the careful estimation of
errors induced by approximating matrix exponentials by first-order
splitting methods. In the present note we show that using higher
order splitting methods allows us to improve the initial estimates
from \cite{OS} to close-to-optimal O$(n^{-(\alpha-\epsilon)})$ bounds for
$\alpha>1$, where $\epsilon>0$ can be chosen arbitrarily small.
\end{abstract}


\section{Introduction}\label{sec1}
The study of approximation rates for Lie-group-valued loops by
polynomial loops is a relatively unexplored topic within the
larger area of nonlinearly constrained approximation. Motivation
is provided by previous density results \cite{La1,La2,PS} for
semi-simple Lie groups, and by more practical needs, e.g., for the
design of para-unitary FIR filters \cite{La1,OMK,Vai}.

In this note, we continue the study of the $\SU(N)$ case, $N\ge
2$, and improve upon the following Jackson-type estimate stated in
\cite{OS}: {\em: For any $\Lip_\alpha$-continuous loop $U:\;\TT
\to \SU(N)$ and $\alpha>1/2$ there exists a sequence of polynomial
loops $U_n:\, \TT\to \SU(N)$ of degree $\le n$ such that the
following asymptotic inequality holds} \be\label{J0}
\|U-U_n\|_C:=\max_{t\in\TT} \|U(t)-U_n(t)\| \le C_{\alpha,N,U}\,
(n+1)^{-\alpha/(1+\alpha)},\qquad n\ge 0. \ee Even though this
estimate is admittedly far from final, to our knowledge it
represented the first nontrivial upper estimate for the achievable
rate of approximation for $\Lip_\alpha$ loops with values in
matrix Lie groups. Note that there is a large gap between the
exponent $\alpha/(\alpha+1)$ established in (\ref{J0}), and the
trivial upper bound $\alpha$ for the maximal order of
approximation of $\Lip_\alpha$ loops following from the classical
Jackson-Bernstein theorems for univariate trigonometric
approximation \cite[Theorem 7.3.3]{DL}.

In the meantime, we realized that some simple modifications in the
proof strategy of \cite{OS} yield close-to-optimal rates, at least
for $\alpha>1$. The major change is to use higher order splitting
methods for exponentials instead of the standard first order approximation
$e^{\sum_{j=1}^J X_j}\approx e^{X_1}\ldots e^{X_J}$.

\begin{theo} Let $N\ge 2$, $\alpha>1$, and $\epsilon>0$. For any
$U(t) \in \Lip_\alpha(\TT\to \SU(N))$, there exists a sequence of
$\SU(N)$-valued polynomial loops $U_n(t)$ of degree $\le n$ such
that \be\label{J1} \|U-U_n\|_C \le C_{\alpha,\epsilon,N,U}\,
(n+1)^{-(\alpha-\epsilon)},\qquad n\ge 0. \ee
\end{theo}
The proof of this result is given below. The restriction
$\alpha>1$ comes from the fact that we currently miss error
formulas for splitting methods of order $k>2$ in terms of higher
order commutators. Another obstacle is the lack of formal proof
for factorizations of $\SU(N)$-valued $\Lip_\alpha$ loops into
exponentials of $\su(N)$-valued $\Lip_\alpha$ loops if $\alpha
\le 1/2$. The latter problem can possibly be removed by using homotopy
arguments as in \cite{La2} (the first author acknowledges
inspiring discussions with W. M. Lawton on this and related
subjects of the present note). See also the remarks at the end of
the next section.

The major open question is whether (\ref{J1}) remains true also with $\epsilon=0$, or if the
nonlinear constraints lead to a slight deterioration of approximation results
compared to the unconstrained case. Settling this question will probably require
a different approach.

\section{Proof of Theorem 1} \label{sec2}
We first recall the facts already proved in \cite{OS}. The
notation we use is either self-explaining or can be found in
\cite{OS} (we have opted to keep notation very close to that of
\cite{OS}, to make the comparison easy). The matrix norm of choice
is the spectral norm.

{\bf Factorization into essentially exponentials of $\su(2)$
loops}. Lemma 4 of \cite{OS} states that for any $\alpha>1/2$ and
any $U(t)\in \Lip_\alpha (\TT\to \SU(N))$, there exist constant
matrices $U_{0,l}\in \SU(N)$ and loops $A_l(t)\in \Lip_\alpha
(\TT\to \su(2))$ such that \be\label{UA} U(t)=\prod_{l=1}^L
U_{0,l} e^{\hat{A}_l(t)},\qquad t\in \TT,\quad L:=N(N-1)/2, \ee
where $\hat{A}_l(t)=T_{ij}A_l(t)$ denotes the canonical extension
of $A_l(t)$ to a $\GG{su}(N)$ loop by the map
$$
A=\left(\ba{cc} a_{11} & a_{12}\\
a_{21}& a_{22} \ea\right) \quad \longmapsto \quad T_{ij}A=
\left(\ba{ccccc} I_{i-1} & 0 & 0 & 0 & 0\\
0 & a_{11} & 0 & a_{21} & 0\\0 & 0 & I_{j-i+1} & 0 & 0 \\
0 & a_{21} & 0 & a_{22} & 0\\0 & 0 & 0 & 0 & I_{N-j} \ea\right)
$$
for some index pair $(i,j)$ with $1\le i< j\le N$ ($I_k$ denotes
the $k\times k$ identity matrix). Moreover, smoothness of the
factors is controlled by smoothness of $U(t)$: \be\label{UALip}
\|A_l\|_{\Lip_\alpha} \le C_{\alpha,N,U}\|U\|_{\Lip_\alpha},\qquad
l=1,\ldots,L. \ee

{\bf Approximation can be done factor-by-factor}. If $P_l(t)$ are
polynomial loops in $\SU(2)$ of degree $\le n$ such that
\be\label{EA} \|e^{A_l(t)}-P_l(t)\|_C \le \epsilon,\qquad
l=1,\ldots,L, \ee for the $A_l(t)$ occurring in the factorization
(\ref{UA}) then
$$
P(t)=\prod_{l=1}^L U_{0,l} \hat{P}_l(t),\qquad
(\;\hat{P}_l(t)=T_{ij}P_l(t)\;)
$$
is a polynomial loop in $\SU(N)$ of degree $\le Ln$ and it
satisfies the estimate \be\label{UAU} \|U(t)-P(t)\|_C \le
L\epsilon. \ee Use Lemma 5 from \cite{OS}.

{\bf Construction of $P_l(t)$}. For any $m>1$, we can approximate
$A_{l}(t)\in \Lip_\alpha(\TT\to \su(2))$ by a
$\su(2)$-valued polynomial loop $R_{l,m}(t)$ of degree $\le m$
at optimal rate (say, by applying Vallee-Poussin means
componentwise), i.e., \be\label{E2} \|
e^{A_l(t)}-e^{R_{l,m}(t)}\|_C \le Cm^{-\alpha},\qquad m>0, \ee see
Lemma 6 a) in \cite{OS}, and
$$
R_{l,m}(t) = \sum_{r=1}^6 \sum_{k=0}^m c_{r,k}B_{r,k}(t)
$$
(to keep notation simple, the dependence of the coefficients
$c_{k,r}$ on $l$ and $m$ is not made explicit, moreover, the terms
with $k=0$ and $r$ even are fake). Here $\{B_{r,k}(t)\}$ is the
designated basis (over $\RR$) for the linear space of dimension
$6m+3$ of all $\su(2)$-valued polynomial loops of degree $\le
m$, see \cite{OS}. Lemmas 1 and 4 in \cite{OS} establish the following facts
which are relevant below:
\be\label{Lip} \|c_{r,k}B_{r,k}(t)\|_C=
|c_{r,k}| \le C(k+1)^{-\alpha},\qquad k=0,\ldots,m, \quad
r=1,\ldots, 6,
\ee
and if ordered properly (e.g.,
lexicographically), products of the form
$$
\prod_{r=1}^6 \prod_{k=0}^m e^{\lambda\, c_{r,k}B_{r,k}(t)}
$$
represent $\SU(2)$-valued polynomial loops of degree $\le 6m$,
independently of the choice of $\lambda>0$ and the set of real
coefficients $\{c_{r,k}\}$. For our applications we take
$\lambda=1/M,$ for some integer $M>1$.

With these preparations, we can write down the final formula for
the approximation of $e^{A_l(t)}$:
\be\label{Pol1} P_l(t):=
\phi(\{c_{r,k}B_{r,k}(t)/M\}_{r=1,\ldots,6;\,k=0,\ldots,m})^M,
\ee
where $\phi(\{X_j\}_{j=1,\ldots,J})$ is a suitable splitting
method for approximating $e^{X_1+\ldots+X_J}$, see the next
paragraph for details. The integers $m$ and $M$ will be fixed
later. Note that the pointwise estimate \be\label{E3}
\|e^{A_l(t)}-P_l(t)\|\le Cm^{-\alpha}+ M \|e^{R_{l,m}(t)/M} -
\phi(\{c_{r,k}B_{r,k}(t)/M\}_{r=1,\ldots,6;\,k=0,\ldots,m})\| \ee
follows from applying the triangle inequality to
$$
e^{A_l(t)}-P_l(t)= (e^{A_l(t)}-e^{R_{l,m}(t)})+ ( (e^{R_{l,m}(t)/M})^M-P_l(t)),
$$
then using (\ref{E2}) for the first term, and Lemma 5 in \cite{OS} for the second. Thus,
the quality of approximation crucially depends on the properties of the chosen splitting method $\phi$.

{\bf Estimate for the second term in (\ref{E3})}. Now we depart from \cite{OS}, where the method $\phi$ of
choice was the first order splitting method
$$
\phi_1(\{X_j\}_{j=1,\ldots,J}):=e^{X_1}e^{X_2}\ldots e^{X_J}.
$$
The error estimate is stated in Lemma 6 b) of \cite{OS}, and leads
to an overall estimate $\le CM^{-1}$ for the second term in the
right-hand side of (\ref{E3}) if $\alpha>1/2$ (and to the
suboptimal asymptotic approximation rate of that paper). We now
show that using higher order symmetric methods leads to
significant improvements. The standard $2$nd order symmetric
method is given by
$$
\phi_2(\{X_j\}_{k=1,\ldots,J}):=e^{X_1/2}\ldots
e^{X_{J-1}/2}e^{X_J}e^{X_{J-1}/2}\ldots e^{X_1/2}.
$$
Following Yoshida  (see \cite{LQ}), symmetric splitting methods of order $2(s+1)$ can be constructed
from a given symmetric method of order $2s$ via the formula
$$
\phi_{2(s+1)}(\{X_j\}_{j=1,\ldots,J}):=
\phi_{2s}(\{a_{s}X_j\}_{j=1,\ldots,J})\phi_{2s}(\{b_sX_j\}_{j=1,\ldots,J})\phi_{2s}(\{a_sX_j\}_{j=1,\ldots,J}),
$$
if one chooses the constants as follows:
$$
a_s=(2-2^{1/(2s+1)})^{-1},\qquad
b_s=-2^{1/(2s+1)}(2-2^{1/(2s+1)})^{-1}.
$$
The order condition for these $\phi_{2s}$ can be stated as follows: For $\lambda\to 0$, we have
$$
\|e^{\lambda(X_1+\ldots + X_J)} -
\phi_{2s}(\{X_j\}_{j=1,\ldots,J})\|=\mathrm{O}(\lambda^{2s+1}).
$$
Using Taylor expansion and rough estimates, the order requirement
translates into the error bound \be\label{Order}
\|e^{\lambda(X_1+\ldots + X_J)} -
\phi_{2s}(\{X_j\}_{j=1,\ldots,J})\|\le C\lambda^{2s+1}
\left(\sum_{j=1}^J \|X_j\|\right)^{2s+1}, \ee valid for
$|\lambda|\le \lambda_{\max}$, with a constant depending on
$\lambda_{\max}$ and $\sum_j \|X_j\|$. More precise error bounds
are available for $\phi_2$, see Remark 2 at the end.

We are now ready to apply this to the family
$\{c_{r,k}B_{r,k}(t)\}_{r=1,\ldots,6;\,k=0,\ldots,m}$ (pointwise
with respect to $t$) with $\lambda=M^{-1}$. By (\ref{Lip}), if
$\alpha>1$ we have \be\label{Sum} \sum_{r=1}^6\sum_{k=0}^m
\|c_{r,k}B_{r,k}(t)\|\le C, \ee where $C$ depends on the
$\Lip_\alpha$-norms of the the $\su(2)$-loops $A_l(t)$, but
not on $m$. Substituting into
(\ref{Order}), we obtain for each $l=1,\ldots,L$ \be\label{E5}
\|e^{R_{l,m}(t)/M} -
\phi_{2s}(\{c_{r,k}B_{r,k}(t)/M\}_{r=1,\ldots,6;\,k=0,\ldots,n})\|_C
\le CM^{-(2s+1)}. \ee From now on, set $\phi=\phi_{2s}$ in the
formula for $P_l(t)$ (and consequently $P(t)$). Substituting into
(\ref{E3}) and taking into account (\ref{EA}), (\ref{UAU}) we
finally arrive at \be\label{E4} \|U(t)-P(t)\|_C\le
L\max_l\,\|e^{A_l(t)}-P_l(t)\|\le C(m^{-\alpha}+M^{-2s}), \ee
where $C$ depends on $s$, $\alpha$, $N$, and on $U(t)$.

{\bf Estimating the degree of $P(t)$}. Consider a large enough
integer $n\ge n_0$ (for $n<n_0$, just use constant $P(t)=I$ to get
the complementing trivial bound $\|U(t)-I\|\le 2$). We will now
fix $m$ and $M$ such that the degree of the above constructed
$P(t)$ is $\le n$ {\it and} the right-hand side in (\ref{E4}) is
asymptotically as small as possible. Due to the recursive
definition of $\phi_{2s}$, the degree of the polynomial loops
$\phi_{2s}(\{c_{r,k}B_{r,k}(t)/M\})$ is bounded by $3^{s-1}$ times
the degree of the polynomial loops
$\phi_{2}(\{c_{r,k}B_{r,k}(t)/M\})$ generated by the $2$nd order
method (for simplicity, we do not indicate the index set
${r=1,\ldots,6;\,k=0,\ldots,m}$
 in the notation). The latter, however, have degree $\le 12m$.
This can be proved as in Lemma 1 of \cite{OS}. Indeed, write
$$
\phi_{2}(\{c_{r,k}B_{r,k}(t)/M\})=
\phi_{1}(\{c_{r,k}B_{r,k}(t)/{2M}\})
\phi_{1}(\{c_{r,k}B_{r,k}(t)^\ast/{2M}\})^\ast.
$$
We already know that the first factor
$\phi_{1}(\{c_{r,k}B_{r,k}(t)/{2M}\})$ has degree $\le 6m$. The
second factor is the Hermitean transpose of
$\phi_{1}(\{c_{r,k}B^\ast_{r,k}(t)/{2M}\})$, and it remains to
check that $\{B^\ast_{r,k}(t)\}$ is such a permutation of the
original basis $\{B_{r,k}(t)\}$ to which Lemma 1 of \cite{OS} can
be applied, leading to the same degree bound.

Putting things together, we see that the degree of $P(t)$ is bounded by $12L 3^{s-1}Mm$. Thus, choosing
the integers $M$, $m$ according to
$$
M = [(12L 3^{s-1})^{-1} n^{\alpha/(\alpha+2s)}], \qquad m = [n^{2s/(\alpha+2s)}],
$$
we guarantee that the degree of $P(t)$ does not exceed $n$. On the other hand,
substitution into (\ref{E4}) yields
$$
 \|U(t)-P(t)\|_C\le Cn^{-2s\alpha/(\alpha+2s)} = Cn^{-\alpha + \alpha^2/(\alpha+2s)}.
$$
This establishes the claim of our theorem, if, for given $\alpha>1$ and $\epsilon>0$, we choose
the order $2s$ of the splitting method large enough.

\medskip
{\bf Remark 1}. There are at least three shortcomings of the
asymptotic estimate (\ref{J1}). First, the constant
$C_{\alpha,\epsilon,N,U}$ depends on $U(t)$ in an unspecified way.
Secondly, the restriction $\alpha > 1$ is mainly due to the use of
the crude error estimate (\ref{Order}) for higher order splitting
methods (see the comments in Remark 2). In addition, for $\alpha
\le 1/2$ the factorization technique of Lemma 4 from \cite{OS}
breaks down (an alternative is addressed in Remark 3). Thirdly, it
is not clear at the moment if one can set $\epsilon=0$ in
(\ref{J1}).

\smallskip
{\bf Remark 2.} For low-order splitting methods such as $\phi_1$
and $\phi_2$, the error bound can be made more precise in terms of
commutator expressions which paves the way for proving (\ref{E5})
(and thus also (\ref{E4})) for some $\alpha\le 1$. For $\phi_1$
this was demonstrated in \cite{OS} (see Lemma 6 for the more
precise error bounds). An improved error bound for the symmetric
$2$nd order method $\phi_2$ has been established in \cite{Su}:
\be\label{Suzuki} \|e^{\lambda(X_1+\ldots + X_J)} -
\phi_{2}(\{X_j\}_{j=1,\ldots,J})\|\le\lambda^3
\Delta(X_{1},\ldots,X_{J}), \ee where
$\Delta(X_1,\ldots,X_J)=\sum_{k=1}^{J-1}\Delta_{2}(X_{k},X_{k+1}+\cdots+X_{J})$,
and
$$\Delta_{2}(A,B)=\frac{1}{12}\{\|[[A,B],B]\|+\frac{1}{2}\|[[A,B],A]\|\}.$$
The advantage is that in our application, where we would apply
(\ref{Suzuki}) to matrix sets of the form $\{c_k
B_{r,k}\}_{k=0,\ldots,m}$ related to the terms of the Fourier
series of a $\su(2)$-valued $\Lip_\alpha$-loop, the sum of the
norms of the appearing triple commutators can be estimated by a
sum of the form $\sum_{k=1}^m (\log k)^2 k^{-3\alpha}$ which
remains uniformly bounded for $\alpha >1/3$ (the details are
worked out in \cite{OS} for the first-order case). In contrast,
using (\ref{Order}) with $s=1$ leads to a constant factor of the
form $(\sum_{k=1}^m k^{-\alpha})^3$.

Unfortunately, generalizations of (\ref{Suzuki}) to Yoshida-type (and any other higher-order)
splitting methods are not known.

\smallskip
{\bf Remark 3}. W. R. Lawton drew our attention to a possible
alternative to loop factorizations of the form (\ref{UA}) proposed
in \cite{OS}. It is well known that $\SU(N)$ is a simply connected
compact $C^\infty$-manifold. Thus, any $\SU(N)$-valued continuous
loop $U(t)$ can be contracted to a point by a homotopy map $\psi:
[0,1]\to C(\TT\to \SU(N))$ (i.e., $\psi$ is continuous,
$\psi(1)=U(t)$, and (without loss of generality) $\psi(0)=I$).
 Let us assume that for $U(t)\in \Lip_\alpha(\TT\to \SU(N))$, the homotopy map $\psi$ can
be found in such a way that $\psi(\xi)\in \Lip_\alpha(\TT\to
\SU(N))$ for all $\xi\in [0,1]$ (i.e., preserves Lipschitz
smoothness along the homotopy path). We do not have a reference
for this assumption but strongly believe that it holds for all
$\alpha>0$.

Now, take a fine enough partition $\xi_0=0<\xi_1<\ldots
<\xi_{K-1}<\xi_K=1$ of $[0,1]$ such that
$\|\psi(\xi_{k-1})-\psi(\xi_k)\|_C\le r_N$, where $r_N$ is the
injectivity radius of the exponential map in the neighborhood of
$I\in \SU(N)$. Then we can write
$$
U(t)=U_1(t)\ldots U_K(t),\qquad U_k(t):=\psi(\xi_{k-1})^\ast \psi(\xi_k),\quad k=1,\ldots,K,
$$
where all $U_k(t)$ belong to $\Lip_\alpha(\TT\to \SU(N))$, and
$$
\|I-U_k\|_C \le  \|\psi(\xi_{k-1})^\ast\|_C \|\psi(\xi_{k-1}-\psi(\xi_k)\|_C=\|\psi(\xi_{k-1})-\psi(\xi_k)\|_C\le r_N,
$$
i.e., $A_k(t) = \log(U_k(t))$ is well defined and belongs to
$\Lip_\alpha(\TT\to \su(N))$ for all $k=1,\ldots,K$. Thus,
\be\label{UA1} U(t) = \prod_{k=1}^K e^{A_k(t)},\qquad
\|A_k\|_{\Lip_\alpha} \le C(\alpha,N,U),\quad k=1,\ldots,K. \ee In
contrast to (\ref{UA}), the number of exponential factors $K$ is
not independent of $\alpha$ and $U$, and the $A_k(t)$ are general
$\su(N)$-valued $\Lip_\alpha$-loops (and not essentially
$\su(2)$-valued as in (\ref{UA})).

However, as long as we accept the dependence on $U(t)$ in the
constant appearing in (\ref{J1}), the factorization (\ref{UA1}) is
sufficient to carry out the above proof with minor changes. The
reduction to the $\SU(2)$ case can be circumvented by working with
a similar basis $\{B_{r,k}(t)\}_{r=1,\ldots,R_N,\,k=0,\ldots,m}$
over $\RR$ for $\SU(N)$-valued polynomial loops of degree $\le m$.
What changes is the number $R_N$ of subsets
$\{B_{r,k}(t)\}_{k=0,\ldots,m}$ of basis elements to be
considered. This number depends only on $N$, and enters the degree
estimates as a linear factor.

\smallskip
{\bf Remark 4}. Loop approximation can be pursued for other Lie
groups and manifolds as well. Work on the closely related case of
$\mathrm{SO}(N)$-valued loops ($N\ge 3$) is ongoing \cite{S}.

\end{document}